\documentclass[12pt,english,reqno]{amsart}
\usepackage[T1]{fontenc}
\usepackage[latin9]{inputenc}
\usepackage{verbatim}
\usepackage{amsthm}
\usepackage{amstext}
\usepackage{amssymb}
\usepackage{cite}
\makeatletter
\numberwithin{equation}{section}
\numberwithin{figure}{section}
\theoremstyle{plain}
\newtheorem{thm}{\protect\theoremname}[section]
  \theoremstyle{definition}
  \newtheorem{defn}[thm]{\protect\definitionname}
  \theoremstyle{definition}
  
  \theoremstyle{plain}
  \newtheorem{lem}[thm]{\protect\lemmaname}
  \theoremstyle{plain}
  \newtheorem{prop}[thm]{\protect\propositionname}
  \theoremstyle{remark}
  \newtheorem{rem}[thm]{\protect\remarkname}
  \theoremstyle{plain}

\usepackage{amscd}
\newtheorem{theorem}{Theorem}[section]
\theoremstyle{plain}

\newtheorem{remark}[theorem]{Remark}

\newcommand{\eqdef}{\stackrel{\mathrm{def}}{=}}

\newcommand{\N}{\mathbb{N}}
\newcommand{\R}{\mathbb{R}}

\newcommand{\dvn}{{\mathrm dv}_{{g}^{(n)}}}

\newcommand{\dv}{{\mathrm dv}_g}

\makeatother

\usepackage{babel}
  \providecommand{\corollaryname}{Corollary}
  \providecommand{\definitionname}{Definition}
  \providecommand{\examplename}{Example}
  \providecommand{\lemmaname}{Lemma}
  \providecommand{\propositionname}{Proposition}
  \providecommand{\remarkname}{Remark}
\providecommand{\theoremname}{Theorem}

\begin{document}

\title[Defect of compactness]{On defect of compactness for Sobolev spaces on manifolds}

\author{Leszek Skrzypczak}

\thanks{One of the authors (L.S.) was supported by National Science Center,
Poland, Grant No. 2014/15/B/ST1/00164.}

\address{Faculty of Mathematics \& Computer Science, Adam Mickiewicz University,
ul. Umultowska 87, 61-614 Pozna\'{n}, Poland}

\email{lskrzyp@amu.edu.pl}

\author{Cyril Tintarev}

\thanks{The other author (C.T.) expresses his gratitude to the Faculty of Mathematics \& Computer Science of Adam Mickiewicz
University for warm hospitality.}

\address{Sankt Olofsgatan 66B, 75330 Uppsala, Sweden}

\email{tammouz@gmail.com}


\maketitle
{\small
\textsc{Abstract.}
\emph{Defect of compactness}, relative to an embedding of two Banach
spaces $E\hookrightarrow F$, is a difference between a
weakly convergent sequence $u_{k}\rightharpoonup u$ in $E$ and $u$, taken up to a remainder that vanishes in the norm of $F$. For Sobolev embeddings in particular, defect of compactness is expressed as a \emph{profile decomposition} - a sum of terms, called \emph{elementary concentrations}, with asymptotically disjoint supports. 
We discuss a profile decomposition for the Sobolev space $H^{1,2}(M)$ of a Riemannian manifold with bounded geometry is a sum of elementary concentrations associated with concentration profiles defined on manifolds different from $M$, that are induced by a limiting procedure. The profiles satisfy an inequality of Plancherel type, and a similar relation, related to the Brezis-Lieb Lemma, holds for $L^p$-norms of profiles on the respective manifolds.} 
\vskip3mm

\section{Introduction}
Defect of compactness of an embedding $E\hookrightarrow F$ of two Banach spaces (the difference between a weakly convergent sequence and its weak limit up to a remainder vanishing in $F$), takes, under general conditions the form of profile decomposition - a sum of, in some sense, decoupled terms, called \emph{elementary concentrations}, which reflect certain asymptotic behavior of the sequence. Profile decomposition for the Sobolev embedding $\dot H^{1,p}(\R^N)\hookrightarrow L^\frac{pN}{N-p}$, $N>p>1$, was found by Solimini \cite{Solimini}, and later, independently, by G\'erard \cite{Gerard} and Jaffard \cite{Jaffard}. It is a sum of decoupled terms of the form $g_kw\eqdef t_k^\frac{N-p}{p}w(t_k(x-y_k))$ with $y_k\in R^N$ and $t_k>0$. \emph{Decoupling} of two rescaling operator sequences, $g_k$ and $g_k'$, of this form refers to $g_k^{-1}g_k'\rightharpoonup 0$, and the asymptotic profile is defined by the inverse rescaling of the original sequence:
$$
g_k^{-1}u_k=t_k^{-\frac{N-p}{p}}u_k(t_k^{-1}x+y_k)\rightharpoonup w.
$$
If the sequence $u_k$ is furthermore bounded in the homogeneous space $H^{1,p}(\R^N)$, it has a profile decomposition with no scaling factors ($t_k=1$) and with remainder vanishing in $L^q$, $p<q<\frac{pN}{N-p}$.
This profile decomposition extends to the Sobolev spaces $H^{1,2}(M)$ of Riemannian manifolds $M$ possessing a rich isometry group, with concentrations taking the form $w\circ \eta_k$, $\eta_k\in \mathrm{Iso}(M)$ and the profiles $w$ defined by $u_k\circ\eta_k^{-1}\rightharpoonup w$ in $H^{1,2}(M)$, see \cite{FiesTin}.  
Profile decomposition for a general embedding $E\hookrightarrow F$  of two Banach spaces, which is cocompact relative to a general group of linear isometries on $E$, is provided in \cite{SoliTi}. 
In relation to Sobolev spaces of Riemannian manifolds, a profile decomposition similar to Solimini's was obtained by Struwe \cite{Struwe}, but only for a specific class of sequences and only on compact manifolds. Elementary concentrations in Struwe's profile decomposition are based on asymptotic profiles defined on the tangent spaces to the manifold at the points of concentration. Applications of Struwe's profile decomposition are elaborated in the monograph \cite{DHR}. Struwe's result was extended to the case of general sequences in a recent paper \cite{deVTi}. 

The result we describe in this announcement generalizes the profile decomposition of \cite{FiesTin} to manifolds that may have no nontrivial isometry group. Complete proofs of all statements here are given in \cite{LSCT}. The problem was proposed to one of the authors several years ago by Richard Schoen \cite{Schoen}.

\section{A ``spotlight'' lemma}
Let $M$  be a smooth, complete $N$-dimensional Riemannian manifold with metric $g$ and a positive injectivity radius $r(M)$. 
We assume that $M$ is
a connected non-compact manifold of bounded geometry.  The latter is defined as follows, e.g. cf. \cite{Shubin}.  
\begin{defn}
\label{def:bg} 
 A smooth Riemannian
manifold $M$ is of bounded geometry if the following two conditions
are satisfied: 

(i) The injectivity radius $r(M)$ of $M$ is positive. 

(ii) Every covariant derivative of the Riemann curvature tensor $R^{M}$of
M is bounded, i.e., $\nabla^{k}R^{M}\in L^{\infty}(M)$ for every
$k=0,1,\dots$ 
\end{defn}
In what follows $B(x,r)$ will denote a geodesic ball in $M$ and
$\Omega_{r}$ will denote the ball in $\R^{N}$ of radius $r$ centered at the origin. Let $r\in(0,r(M))$ be fixed. Then the Riemannian exponential map $\mathrm{exp}_{x}$ is a diffeomorphism of $\{v\in T_xM:\, g_x(v,v)< r\}$ onto $B(x,r)$.  For each $x\in M$ we choose an orthonormal basis for $T_xM$ which yields an identification $i_x:\R^N \rightarrow T_xM$.  Then  $e_{x}:\Omega_{r}\to B(x,r)$ will denote  a geodesic normal coordinates at $x$ given by $e_{x}=\mathrm{exp}_{x}\circ i_x$. We do not require smoothness of the map $i_x$ with respect to $x$. We will consider the maps $e_x$ as defined on the balls $\Omega_a$ with $a=\frac34 r(M)$.

\begin{defn}
\label{def:discr}  A subset $Y$ of  Riemannian manifold $M$ is called  $\varepsilon$-discretization of $M$, $\varepsilon>0$,  if the distance between any two distinct points of $Y$ is greater than or equal to $\varepsilon$ and 
\[ M = \bigcup_{y\in Y} B(y,\varepsilon).\]
\end{defn}
Any connected Riemannian manifold $M$ has a $\varepsilon$-discretization for any $\varepsilon>0$, and if $M$ is of bounded geometry then for any $t\ge  1$ the covering  $\{ B(y,t\varepsilon)\}_{y\in Y}$ is uniformly locally finite.   

The following statement is a counterpart of the "cocompactness" lemma proved in the Euclidean case by Lieb \cite{Lieb}.
\begin{lem}["Spotlight lemma"]
\label{lem:spotlight}Let $M$ be an $N$-dimensional Riemannian manifold of bounded geometry and let $Y\subset M$ be a $r$-discretization of $M$,  $r<r(M)$. 
 Let  $(u_{k})$ be a bounded sequence in $H^{1,2}(M)$.  Then, $u_{k} \to 0$ in $L^{p}(M)$ for any $p\in(2,\frac{2N}{N-2})$ if and only if  
 $u_{k}\circ e_{y_{k}}\rightharpoonup 0$ in
$H^{1,2}(\Omega_{r})$  for any  sequence $(y_k)$, $y_{k}\in Y$.  
\end{lem}
\section{Manifold-at-infinity}
In what follows we consider the radius $\rho<\frac{r(M)}{8}$ and a $\hat{\rho}$-discretization  $Y$ of $M$, $\frac{\rho}{2}<\hat{\rho}<\rho$. We will write $\N_{0}\eqdef\N\cup\lbrace0\rbrace$.
\begin{defn}
\label{def:trailing} Let $(y_{k})_{k\in\N}$ be a sequence in $Y$ that is an enumeration of the infinite subset of $Y$.
A countable family $\lbrace(y_{k;i})_{k\in\N}\rbrace_{i\in\N_{0}}$
of sequences on $Y$ is called  \emph{a  trailing system} for $(y_{k})_{k\in\N}$ if for every $k\in\N$ $(y_{k;i})_{i\in\N_{0}}$
is an ordering of $Y$ by the distance from $y_{k}$, that is, an
enumeration of $Y$ such that $d(y_{k;i},y_k)\le d(y_{k;i+1},y_k)$ for
all $i\in\N_{0}$. In particular, $y_{k;0}=y_k$. 
\end{defn}
It is easy to see that any enumeration of the infinite subset of  $Y$ admits a trailing system.
With a given trailing system $\lbrace(y_{k;i})_{k\in\N}\rbrace_{i\in\N_{0}}$ we associate a manifold $M_\infty^{(y_{k;i})}$ defined by a gluing data described below.  

We give only a rough sketch of construction that involves many technical details. 
For each such  pair $(i,j)\in\N_0$ we consider maps on $\bar{\Omega}_{2\rho}$
\[
\psi_{ij,k}\eqdef e_{y_{k;i}}^{-1}\circ e_{y_{k;j}},
\]
The range of $e_{y_{k;j}}$, which is $\bar B(y_{k;j},2\rho)$, may not necessarily fall into the domain of $e_{y_{k;i}}^{-1}$, which is 
$B(y_{k;i},a)$, so the maps $\psi_{ij,k}$ are defined only for a subset of $i$, $j$ and $k$. There is, however, a certain non-empty set $\mathcal K\subset \N_0\times \N_0$ such that  $\psi_{ij,k}\eqdef e_{y_{k;i}}^{-1}\circ e_{y_{k;j}}$ is a map $\bar{\Omega}_{2\rho}\to\Omega_{a}$ for all $k$ sufficiently large whenever  $(i,j)\in\mathcal K$. 

From boundedness of the geometry and the Ascoli-Arzela theorem, it follows that there is a renamed subsequence of $(\psi_{ij,k})_{k\in\N}$, $(i,j)\in\mathbb K$, that converges
in $C^{\infty}(\bar{\Omega}_{2\rho})$ to some smooth function $\psi_{ij}:\bar{\Omega}_{2\rho}\to\Omega_{a}$,
and, moreover, we may assume that the same extraction of $(\psi_{ji,k})_{k\in\N}$
converges in $C^{\infty}(\bar{\Omega}_{2\rho})$ as well. 
We define $\Omega_{ij}\eqdef\psi_{ij}(\Omega_{\rho})\cap\Omega_{\rho}$.
This set may generally be empty. Let us define 
\begin{equation}
\mathbb{K}\eqdef \{(i,j)\in \mathcal K:\;\Omega_{ij}\neq\emptyset\}.
\end{equation}
 Basing on the gluing theorem in \cite{Gallier3} we can associate the family of sets 
 $\Omega_{ij}$ and the maps $\psi_{ij}$ with a differentiable manifold as follows.
 
 \begin{prop}
Let $M$ be a Riemannian manifold with bounded geometry and let $Y$ be its discretization. 

For any trailing system $\lbrace(y_{k;i})_{k\in\N}\rbrace_{i\in\N_{0}}$ related to the sequence $(y_k)$  in $Y$ there exists a smooth manifold $M_{\infty}^{(y_{k;i})}$  with an atlas $\{(U_i,\tau_i)\}_{i\in\N_{0}}$ such that:\\
1)   $\tau_i(U_i)=\Omega_\rho$, \\
2)  there exists a renamed subsequence of $k$ such that for any two charts  $(U_i,\tau_i)$ and $(U_i,\tau_i)$ with $U_i\cap U_j\not=\emptyset$  the corresponding  transition map $\psi_{ij}: \tau_j(U_j\cap U_i)\rightarrow \tau_i(U_j\cap U_i)$ is the following $C^\infty$-limit: 
\[
\psi_{ij}=\lim_{k\rightarrow \infty} e_{y_{k;i}}^{-1}\circ e_{y_{k;j}}. 
\]
\end{prop}
For convenience we introduce "inverse" charts $\varphi_i= \tau_i^{-1}$ so that $\varphi_{j}^{-1}\circ\varphi_{i}=\psi_{ji}:\Omega_{ij}\to\Omega_{ji}$.

We can endow the manifold $M_{\infty}^{(y_{k;i})}$ with a metric which is  also related to the asymptotic properties of $M$. For any $i\in \N_0$ we define a metric tensor $g^{(i)}$ on $\Omega_\rho$, by the limiting procedure on a suitable renamed subsequence:
\begin{equation}
\widetilde{g}_{\xi}^{(i)}(v,w)\eqdef\lim_{k\to\infty} g_{e_{y_{k;i}(\xi)}}   \big(d e_{y_{k;i}}(v), d e_{y_{k;i}}(w)\big) , \; \xi\in \Omega_\rho \; \text{and}\; v,w\in \R^N,\label{eq:pre-metric}
\end{equation}
Existence of the limit follows from the boundedness of the geometry of  $M$.
Afterwards we pull the metric tensor back onto $U_i=\varphi_i(\Omega_\rho)\subset M_{\infty}^{(y_{k;i})}$ via $\varphi_i^{-1}$: 
\begin{align} \label{eq:metric}
\widetilde{g}_{x}(v,w)  &\eqdef \widetilde{g}^{(i)}_{\varphi_i^{-1}(x)}\big(d \varphi_i^{-1}(v), d \varphi_i^{-1}(w)\big) , \\ 
\qquad &  x\in \varphi_{i}(\Omega_\rho) \subset M_{\infty}^{(y_{k;i})}  \; \text{and}\; v,w\in T_x M_{\infty}^{(y_{k;i})}, \nonumber 
\end{align}
and then verify compatibility of the definition on overlapping charts.   
\begin{defn}
A manifold at infinity $M_{\infty}^{(y_{k;i})}$ of a manifold $M$
with bounded geometry, generated by a trailing system $\lbrace(y_{k;i})_{k\in\N}\rbrace_{i\in\N_{0}}$ of a sequence $(y_k)$ in $Y$,
is a differentiable manifold 
supplied with a Riemannian metric tensor $\widetilde{g}$ defined  by (\ref{eq:metric}).   
\end{defn}
Since all the limits in construction are uniform $C^\infty$-limits, manifolds at infinity of $M$ are also of bounded geometry.
\begin{prop}
Let $M$ be a Riemannian  manifold with bounded geometry and let $Y$ be its $\hat\rho$-discretization, $\rho/2<\hat{\rho}<\rho < \frac{r(M)}{8}$. Then for every discrete sequence
$(y_{k})$ in $Y$ and its trailing system $\lbrace(y_{k;i})_{k\in\N}\rbrace_{i\in\N_{0}}$ there exists a renamed subsequence $(y_{k})$  that generates a Riemannian manifold at infinity $M_{\infty}^{(y_{k;i})}$ of the manifold $M$. The manifold $M_{\infty}^{(y_{k;i})}$ has bounded geometry with an injectivity radius not less than $\rho$. 
\end{prop}

\begin{rem}\label{rem:bg}
%
If $M'$ is another manifold such that $M$ and $M'$ have respective
compact subset $M_{0}$ and $M_{0}'$ such that $M\setminus M_{0}$
is isometric to $M'\setminus M_{0}'$, i. e. if $M'$ is $M$ up to
a compact perturbation, then  their manifolds at infinity for corresponding trailing systems coincide.
From this follows that  manifold at infinity of the manifold $M$  is not necessarily diffeomorphic 
to $M$. 

\end{rem}

\section{Concentration profiles. The main result}
\begin{defn}
Let $M$ be a manifold of bounded geometry and $Y$ be its discretization. Let $(u_{k})$ be a bounded
sequence in $H^{1,2}(M).$ Let $(y_{k})$ be a sequence of points
in $Y$
and let $\lbrace(y_{k;i})_{k\in\N}\rbrace_{i\in\N_{0}}$
be its trailing system. 
One says that $w_{i}\in H^{1,2}(\Omega_{\rho})$ is  \emph{a local profile}
of $(u_{k})$ relative to a trailing sequence $(y_{k;i})_{k\in\N}$, if,
on a renamed subsequence, $u_{k}\circ e_{y_{k;i}}\rightharpoonup w_{i}$
 in $H^{1,2}(\Omega_{\rho})$ as $k\to\infty$. If $(y_{k})$ is a
renamed (diagonal) subsequence such that $u_{k}\circ e_{y_{k;i}}\rightharpoonup w_{i}$
in $H^{1,2}(\Omega_{\rho})$ as $k\to\infty$ for all $i\in\N_{0}$, then
the family $\lbrace w_{i}\rbrace_{i\in\N_{0}}$
is called an \emph{array of local profiles of $(u_{k})$ relative to the trailing system $\lbrace(y_{k;i})_{k\in\N}\rbrace_{i\in\N_{0}}$ of the sequence $(y_{k})$}.
\end{defn}

\begin{prop}
\label{prop:W}Let $M$ be a manifold of bounded geometry and let $Y$ its discretization. Let 
$(u_{k})$ be a bounded sequence in $H^{1,2}(M).$  Let $\lbrace w_{i}\rbrace_{i\in\N_{0}}$
be an array of local profiles of $(u_{k})$ associated with a trailing system $\lbrace(y_{k;i})_{k\in\N}\rbrace_{i\in\N_{0}}$ related to the sequence $(y_k)$ in
$Y$. Then there exists a function $w:\ M_{\infty}^{(y_{k;i})}\to\R$
such that $w\circ\varphi_{i}=w_{i}$, $i\in\N_{0}$, where $\varphi_{i}:\Omega_{\rho}\to M_{\infty}^{(y_{k;i})}$ are local coordinate maps of 
$ M_{\infty}^{(y_{k;i})}$.
\end{prop}

\begin{defn}\label{def:gp}
Let $\lbrace w_{i}\rbrace_{i\in\N_{0}}$ be a local profile array of a bounded sequence $(u_{k})$
in $H^{1,2}(M)$ relative to a trailing system  $\lbrace(y_{k;i})_{k\in\N}\rbrace_{i\in\N_{0}}$. The function
$w:M_{\infty}^{(y_{k;i})}\to\R$ given by Proposition \ref{prop:W}
is called the \emph{global profile} of the sequence $(u_{k})$ relative
to $(y_{k;i})$.
\end{defn}

Since $M$ has bounded geometry, we may fix a uniformly smooth partition of unity $\lbrace\chi_{y}\rbrace_{y\in Y}$ subordinated to the uniformly finite covering of $M$ by geodesic
balls $\lbrace B(y,\rho)\rbrace_{y\in Y}$. 

\begin{defn}
Let $M$ be a manifold of bounded geometry and let $Y$ be its discretization. Let $M_{\infty}^{(y_{k;i})}$   be a manifold at infinity of $M$ generated by a trailing system  $\lbrace(y_{k;i})_{k\in\N}\rbrace_{i\in\N_{0}}$. 
An \emph{elementary concentration} associated with a function $w:M_{\infty}^{(y_{k;i})}\to\R$ is a sequence $(W_{k})_{k\in\N}$ of functions $M\to\R$ given
by
\begin{equation}
W_{k}=\sum_{i\in\N_{0}}\chi_{y_{k;i}}w\circ\varphi_{i}\circ e_{y_{k;i}}^{-1},\qquad k\in\N.\label{eq:elem-conc}
\end{equation}
where   $\varphi_{i}$ are the local coordinate maps of manifold $M_{\infty}^{(y_{k;i})}$.
\end{defn}

In heuristic terms, after we find limits $w_{i}$, $i\in\N_{0},$ of  the sequence $(u_{k})$ under the ``trailing spotlights'' $(e_{y_{k;i}})_{k\in\N_{0}}$
that follow different trailing sequences $(y_{k;i})_{k\in\N}$ of $(y_{k})$, we give an approximate reconstruction $W_{k}$ of $u_{k}$
``centered'' on the moving center $y_{k}$ of the ``core spotlight''.
We do that by first splitting $w$ into local profiles $w\circ\varphi_{i}$, $i\in\N_0$, on the set $\Omega_\rho$, casting them onto the manifold $M$ in the vicinity of $y_{k;i}$ by composition with $e_{y_{k;i}}^{-1}$,
and patching all such compositions together by the partition of unity
on $M$. Such reconstruction approximates $u_{k}$ on geodesic balls
$B(y_{k},R)$ with any $R>0$, but it ignores the values of $u_{k}$
for $k$ large on the balls $B(y_{k}',R)$, with $d(y_{k},y_{k}')\to\infty$,
where $u_{k}$ is approximated by a different local concentration.
It has been shown in \cite{FiesTin} for the case of manifold $M$ with  cocompact action of a group of isometries
(in particular, for homogeneous spaces) that a global reconstruction
of $u_{k}$, up to a remainder vanishing in $L^p(M)$, is a sum elementary concentrations associated with all such mutually decoupled sequences. 

Similarly, the profile decomposition theorem below, which is the main result
of this paper, says that any bounded sequence $(u_{k})$ in $H^{1,2}(M)$
has a subsequence that, up to a remainder vanishing in $L^{p}(M)$, $p\in(2,2^*)$, equals a sum of decoupled elementary concentrations. 
 To simplify the notation we will index the  sequences, the related trailing systems  and the corresponding manifold by $n$, i.e. below we  write     $y^{(n)}$, $y^{(n)}_{k;i}$, and  $M^{(n)}_\infty$.      

\begin{thm}
\label{thm:main}
Let $M$ be a manifold of bounded geometry and let $Y$ be its discretization.  Let
$(u_{k})$ be a sequence in $H^{1,2}(M)$ weakly convergent to some function $w^{(0)}$ in $H^{1,2}(M)$. 
Then there exists a renamed subsequence of $(u_{k})$, sequences $(y_{k}^{(n)})_{k\in\N}$ in $Y$ 
, and associated with them global profiles $w^{(n)}$ 
on the respective manifolds at infinity  
$M_{\infty}^{(n)}$, 
$n\in\N$, such that $d(y_{k}^{(n)},y_{k}^{(m)})\to\infty$ when $n\neq m$, and 
\begin{equation}
u_{k}-w^{(0)}-\sum_{n\in\N}W_{k}^{(n)}\to 0\mbox{ in }L^{p}(M),\; p\in(2,2^{*}),\label{eq:PD}
\end{equation}
where $W_{k}^{(n)}=\sum_{i\in\N_{0}}\chi_{i}^{(n)}w^{(n)}\circ\varphi_{i}^{(n)}\circ e_{y_{k;i}^{(n)}}^{-1}$ are elementary concentrations,  $\varphi_{i}^{(n)}$ are the local coordinates of the manifolds $M^{(n)}_\infty$ and $\{\chi_i^{(n)}\}_{i\in\N_0}$ the corresponding partitions of unity. 
The series $\sum_{n\in\N}W_{k}^{(n)}$ converges in $H^{1,2}(M)$
unconditionally and uniformly in $k\in\N$. Moreover,
\begin{equation}
\|w^{(0)}\|_{H^{1,2}(M)}^{2} +\sum_{n=1}^{\infty}\|w^{(n)}\|_{H^{1,2}(M^{(n)}_\infty)}^{2}\le  \limsup\|u_{k}\|_{H^{1,2}(M)}^{2}\ ,\label{eq:Plancherel}
\end{equation}
 and 
\begin{equation}
\int_{M}|u_{k}|^{p}d\dv\to\int_{M}|w^{(0)}|^{p}\dv+\sum_{n=1}^{\infty}\int_{M_{\infty}^{(n)}}|w^{(n)}|^{p}\dvn.
\label{eq:newBL}
\end{equation}
\end{thm}

\begin{remark} In \cite{LSCT} it is shown that Theorem~\ref{thm:main}
implies the profile decomposition of \cite{FiesTin} in the case when $M$ is cocompact relative to a discrete isometry group.
\end{remark}
It is interesting to compare the objects at infinity in the profile decomposition \eqref{eq:PD} and in the profile decompositions in \cite{Struwe,deVTi}. In the latter, loss of compactness occurs due to  blowup concentrations, and concentration profiles, defined by behavior of the sequence near a given point, are functions on the tangent space, which can be seen as the manifold-at-infinity created by the concentration mechanism at work - zooming into the manifold $M$ at a given point. In \eqref{eq:PD} concentration profiles are generated by localized shifts to infinity, followed by a reassembly on a new manifold.  
In both profile decompositions, sum of the energies of profiles on respective manifolds at infinity is dominated by the energy of the sequence, and analogous relations hold for the $L^p$-norms.


\begin{thebibliography}{10}




\bibitem{deVTi} G. Devillanova, C. Tintarev, On defect of compactness on compact Riemannian manifolds, in preparation

\bibitem{DHR}O. Druet, E. Hebey, F. Robert, Blow-up Theory for Elliptic
PDEs in Riemannian Geometry, Princeton University Press, 2004.


\bibitem{FiesTin}K.-H. Fieseler, K. Tintarev, Semilinear elliptic
problems and concentration compactness on non-compact Riemannian manifold,
J. Geom, Anal. \textbf{13}, (2003), 67--75

\bibitem{Gallier}J. Gallier, J. Quaintance, Notes on Differential
Geometry and Lie Groups, Book in progress (2017), http://www.cis.upenn.edu/\textasciitilde{}jean/gbooks/manif.html

\bibitem{Gallier3} J. Gallier, D. Xu, M. Siqueira, Parametric pseudo-manifolds,
Differential Geometry and its Applications 3\textbf{0} (2012) 702--736.

\bibitem{Gerard} P. G\'erard, Descriptio de compacit\'e de l'injection de
Sobolev. ESAIM Control Optim. Calc. Var. \textbf{3} (1998 )213--233.

\bibitem{Jaffard} S. Jaffard, Analysis of the lack of compactness
in the critical Sobolev embeddings, J. Funct. Analysis. \textbf{161} (1999), 384--396. 

\bibitem{Lieb}E. Lieb, On the lowest eigenvalue of the Laplacian
for the intersection of two domains. Invent. Math. \textbf{74} (1983), 441--448. 
\bibitem{Schoen} Richard Schoen: personal communication, Stanford, winter 2003.

\bibitem{Shubin} M. A. Shubin, Spectral theory of elliptic operators
on noncompact manifolds, Méthodes semi-classiques, Vol. 1 (Nantes,
1991), Astérisque 207, 35--108 (1992).

\bibitem{LS} L. Skrzypczak, Atomic decompositions on manifolds with bounded geometry, Forum Math. \textbf{10}(1998), 19--38.

\bibitem{LSCT}L. Skrzypczak, C. Tintarev, Defect of compactness for Sobolev spaces on manifolds with bounded geometry, preprint. 

\bibitem{Solimini} S. Solimini, A note on compactness-type properties
with respect to Lorentz norms of bounded subsets of a Sobolev space,
Ann. Inst. H. Poincar\'e Anal. Non-lin\'eaire
\textbf{12} (1995), 319--337. 

\bibitem{SoliTi} S.~Solimini, C. Tintarev, Analysis of concentration
in the Banach space, Comm. Contemp. Math. \textbf{18} (2016), 1550038 (33 pages).


\bibitem{Struwe}M. Struwe, A global compactness result for elliptic
boundary value problems involving limiting nonlinearities, Math. Z.
\textbf{187}, 511--517 (1984).
 

\end{thebibliography}
\end{document}